\documentclass[12pt]{article}
\usepackage{array}
\usepackage{amsfonts,amssymb,amsmath,amsthm}
\usepackage{enumerate}
\usepackage[all]{xy}
\usepackage{hyperref}
\newtheorem{theorem}{Theorem}[section]
\newtheorem{proposition}[theorem]{Proposition}
\newtheorem{proposition-definition}[theorem]{Proposition-Definition}

\newtheorem{corollary}[theorem]{Corollary}
\theoremstyle{remark}

\theoremstyle{definition}

\DeclareMathOperator{\Hom}{Hom}

\newcommand{\Mod}{\mbox{-}\mathrm{Mod}}
\newcommand{\N}{\mathbb{N}}
\newcommand{\Q}{\mathbb{Q}}

\DeclareMathOperator{\Tor}{Tor}

\newcommand{\Z}{\mathbb{Z}}

\newcommand{\bcat}{\mathcal{B}}

\newcommand{\cat}{\mathcal{C}}

\newcommand{\codiag}{\nabla}
\DeclareMathOperator{\coker}{\mathrm{Coker}}

\newcommand{\egamma}{e}

\newcommand{\ggr}{\mbox{-}\Gamma\mbox{-gr}}

\newcommand{\id}{\mathrm{id}}
\newcommand{\image}{\operatorname{Im}}
\newcommand{\im}{\operatorname{Im}}

\newcommand{\md}{\mbox{-Mod}}
\newcommand{\ob}{\mathrm{Ob}}
\newcommand{\onto}{\twoheadrightarrow}

\DeclareMathOperator{\rad}{\operatorname{Rad}}
\newcommand{\smallsub}{\Subset}
\DeclareMathOperator{\supp}{supp}

\author{Ana Paula Santana\thanks{\scriptsize University of Coimbra, CMUC, Department of
Mathematics, 3001-501 Coimbra, Portugal\newline aps@mat.uc.pt}\phantom{,}\ and Ivan Yudin\thanks{\scriptsize University of Coimbra, CMUC, Department of
Mathematics, 3001-501 Coimbra, Portugal\newline yudin@mat.uc.pt}}
\title{
On functors preserving projective resolutions
}
\date{}

\begin{document}

\maketitle
\begin{abstract}
It is important for applications of Homological Algebra in Representation
Theory to have control over the behaviour of (minimal) projective resolutions
under various functors.
In this article we describe three broad families of functors that preserve such
resolutions. We will use these results in our  work on Representation Theory
of Schur algebras.
\medskip

\medskip

\end{abstract}
\section*{Introduction}

It is well known that finite dimensional modules over finite dimensional
algebras admit minimal projective resolutions. The explicit knowledge of  a minimal
projective resolution
for a given module permits to reduce computation of Ext-groups involving this
module to an, albeit sometimes complicated, linear algebra problem.

The present article was conceived in the sequence of our research of
 Ext-groups between Weyl modules for the general linear group $GL_n$ or, equivalently, of Ext-groups between Weyl modules for Schur algebras. The determination of
Ext$^i_{GL_n}(M,N)$, where $M$, $N$ are Weyl modules, is an open problem in the representation theory
of $GL_n$  and of other algebraic groups, see \cite{V}, \cite{P} for the cases of $SL_2$ and $GL_2$.
We hope that the construction of suitable projective resolutions of Weyl modules over Schur algebras can shed some light on this problem.
%In \cite{P}  A. Parker determined  recursive formulas  which  relate the cohomology of certain Weyl modules to that of Weyl modules having smaller highest weight and allow the computation of these Ext-groups for $SL_2$. Following this, K. Erdmann, K. Hannabuss, A. Parker in \cite{K} derive various generating functions for these groups. V. Miemietz and W. Turner in \cite{V}

In our work to explicitly  build  these  projective resolutions, we found that it is important to be able to pass
between different abelian categories with suitably chosen functors. These
functors  should preserve minimal projective resolutions, not
necessarily for all objects,  but at least for the objects we are interested in.

In this article we present our results concerning three of these functors.
Let $\Gamma$ be a monoid.
The first two functors are  defined on the category
of $\Gamma$-graded modules over a $\Gamma$-graded algebra.
The first functor is the forgetful functor that erases the grading information.
The second functor is defined as a twisted product $- \ltimes_{\Gamma} N$ for a
$B$-module $N$, where $B$ is a $\Gamma$-algebra (the explicit definitions are given in Section~3).
In Sections 2 and  3 we determine sufficient conditions for which these
functors preserve minimal projective resolutions for all objects.

In the last section we study the functor $A/I \otimes_A -$, where $A$ is an
algebra and $I$ is an ideal. This functor usually does not preserve (minimal)
projective resolutions for all objects, but in favorable circumstances it
preserves (minimal) projective resolutions of $A/I$-modules considered as
$A$-modules. In this case we say that $I$ is a \emph{stratifying} ideal. The
equivalent condition for $I$ to be stratifying is that $\Tor^A_k(A/I,A/I)=0$ for
$k\ge 1$.

If $R$ is a subring of $A$ which does not necessarily lie inside the center of
$A$, one can define relative torsion groups $\Tor_k^{(A,R)}(A/I,A/I)$
(see~\cite{homology}). In~\cite{LU} we found a combinatorial criterion for
$\Tor_k^{(A,R)}(A/I,A/I)=0$, $k\ge 1$.
In the last section of the article we establish a sufficient condition on the
triple $(A,R,I)$ under which $\Tor^{(A,R)}_k(A/I,A/I)=0$ for $k\ge 1$ implies
the same for $\Tor^A_k(A/I,A/I)$. Thanks to this, one can use the combinatorial
criterion from~\cite{LU} to prove that $I$ is a stratifying ideal in $A$.

For the convenience of the reader, in the first section of the article, we prove some results concerning minimal
projective resolutions in general abelian categories. Despite the fact that these results are
sometimes accepted as true in this general context, we were not able to find
explicit proofs for them in the literature. So we chose to provide full
details in their treatment.

\section{Superfluous subobjects and minimal projective covers in abelian
categories}
In this section we work in a very abstract setting, collecting definitions and proving  results in general abelian categories.
These results and concepts will be then applied to concrete categories in the following sections of the article.

Let $\cat$ be an abelian category.
We follow MacLane~\cite{cwm}, and define
 a \emph{subobject} of $Y \in \cat$ as an equivalence class of
monomorphisms $\psi \colon X \to Y$, where $\psi \sim \psi' \colon X'
\to Y$ if there is an isomorphism $\rho \colon X \to X' $ such that $\psi'\circ
\rho = \psi$.
We use
upper case letters to denote objects and italic uppercase letters to denote subobjects. To indicate
that
$\mathcal{X}$ is a subobject of $Y$ we write $\mathcal{X} \subset Y$.
If $\mathcal{X}$ and $\mathcal{Z}$ are subobjects of $Y$ represented by
monomorphisms $\psi_X \colon X \to Y$ and $\psi_Z \colon Z \to Y$ we say that
$\mathcal{Z}$ is contained in $\mathcal{X}$ and write $\mathcal{Z} \le
\mathcal{X}$ if there is $\phi \colon Z \to X$ such that $\psi_X \circ \phi =
\psi_Z$.
Given an object $Y$ we denote by $\mathcal{Y}$ the maximal subobject of
$Y$, i.e. the subobject given by the equivalence class of $\id_Y$.

Dually a \emph{quotient} of $Y$ is an equivalence class of epimorphisms
$\tau\colon Y \to Z$, where $\tau \sim \tau' \colon Y\to Z'$ if and only if
there is an isomorphism $\rho \colon Z \to Z'$ such that $\rho \circ \tau =
\tau'$.

We adopt the convention that for an arrow $\phi\colon X \to Y$ the kernel
$\ker \phi$ and the image $\im \phi$ of $\phi$ are subobjects of $X$ and
$Y$, respectively. Similarly
the cokernel $\coker \phi$ of $\phi$ is a quotient of $Y$. This is a legitimate
point of view as explained in~\cite[VIII.1]{cwm}.

As usual  we use $i_k$ to denote the canonical embeddings associated with
direct sums. For an object $A$ we denote by  $\codiag_A$ the  codiagonal map, which is
determined by the universal property of the direct sum
\begin{equation*}
\begin{gathered}
\xymatrix{
0 \ar[r] \ar[d] & A \ar[d]^{i_1} \ar@/^4ex/[rdd]^{\id_A} \\
A \ar[r]^{i_2} \ar@/_3ex/[rrd]_{\id_A} & A\oplus A \ar@{.>}[rd]^{\exists!\,\codiag_A} \\
&& A .
}
\end{gathered}
\end{equation*}
Given two subobjects $\mathcal{X}_1$ and $\mathcal{X}_2$ of $Y \in \cat$, we
define the  sum $\mathcal{X}_1 + \mathcal{X}_2$  as the image of $\codiag_Y \circ ( \psi_1 \oplus
\psi_2)$ in $Y$,
where, for $i=1$, $2$, the map  $\psi_i$ is a  representative of
$\mathcal{X}_i$.

For a morphism $\phi \colon X \to Y$ and a subobject $\mathcal{N}$ of $X$
represented by a monomorphism $\psi_N\colon N \to X$, we define the subobject
$\phi(\mathcal{N})$ to be  $\im (\phi \circ \psi_N)$.

Now, if $\mathcal{M}$ is a subobject of $Y$ represented by a monomorphism
$\psi_M\colon M \to Y$, we define $\phi^{-1} (\mathcal{M}) \subset X$
as the image of the left vertical arrow in the pull-back diagram
\begin{equation}
\begin{gathered}
\xymatrix{
\widetilde{M} \ar[r] \ar[d] & M \ar[d]^{\psi_M} \\
X \ar[r]^\phi & Y .
}
\end{gathered}
\end{equation}

 Suppose $\mathcal{X}\subset Y$. We say that $\mathcal{X}$ is a \emph{superfluous subobject} of $Y$ if for
any subobject $\mathcal{T}$ of $Y$ the equality $\mathcal{X}+\mathcal{T}
=\mathcal{Y}$ implies $\mathcal{T}=\mathcal{Y}$.
We will write $\mathcal{X}\smallsub Y$ in this case.

We will define minimal projective covers using superfluous subobjects. As in
this article we study preservability of projective covers under various
functors, it is convenient to have at hand various elementary properties of
superfluous subobjects, which will be given in Proposition~\ref{longlist}. These
properties are well-known and easy to prove in the case $\cat$ is the category of
modules over a ring. To lift these properties to arbitrary abelian
categories, we will use the Freyd-Mitchell embedding:

\begin{theorem}\label{freyd}
Let $\bcat$ be a small abelian category. Then there exists a ring $R$ such that
there is a full and faithful exact functor $F\colon \bcat \to R\Mod$.
\end{theorem}
Notice that the above result refers only to small abelian categories. As we work
with arbitrary abelian categories, we will use the following fact proved on
page 85 of~\cite{ca2}.
\begin{theorem}\label{borceux}
Let $\cat$ be an abelian category and $X$ a set of objects in $\cat$. Then
there is a full subcategory $\bcat(X)$ of $\cat$ such that
\begin{enumerate}[i)]
\item $X \subset \ob \bcat(X)$;
\item $\bcat(X)$ is  a small abelian category;
\item $\bcat(X)$ is stable under finite limits and colimits.
\end{enumerate}
\end{theorem}

We say that a functor  $F$ between abelian categories  is exact if  one of the two
equivalent conditions holds:
\begin{enumerate}[i)]
\item $F$ preserves exact sequences;
\item
 $F$ preserves
finite limits and finite colimits.
\end{enumerate}

Below we list standard properties of fully faithful exact functors, which we
will use without further reference.
If $F\colon \bcat\to \cat$ is a full and faithful  exact
functor between abelian categories, then
\begin{enumerate}[$\bullet$]
\item $F(0)\cong 0$;
% \item $F$ reflects isomorphisms;
\item $F$ preserves and reflects  monomorphisms and epimorphisms;
% \item $F$ reflects monomorphisms and epimorphisms;
\item $F$ preserves and reflects kernels, cokernels, and, thus, also images.
\end{enumerate}
Given $\mathcal{X} \subset Y$ in $\bcat$, we define $F(\mathcal{X}) \subset
F(Y)$ as the equivalence class of $F(\psi_X)$, where $\psi_X \colon X \to Y$ is a
representative of $\mathcal{X}$. Then there hold
\begin{enumerate}[$\bullet$]
\item If $\mathcal{X}_1$, $\mathcal{X}_2 \subset Y$, then $F(\mathcal{X}_1 +
\mathcal{X}_2) = F(\mathcal{X}_1) + F(\mathcal{X}_2)$;
\item If $\mathcal{X}_1$ and  $\mathcal{X}_2$ are two different subobjects in
$Y$ then $F(\mathcal{X}_1) \not=F(\mathcal{X}_2)$;
\item  $F$ preserves images and pre-images of subobjects.
\end{enumerate}
We can now state and prove the properties of superfluous subobjects that
we mentioned before.
\begin{proposition}
\label{longlist}
Let $\cat$ be an abelian category.
Then the following statements hold:
\begin{enumerate}[(i)]
\item \label{small_homomorphism}
		Suppose that $\mathcal{N}\smallsub M$  and $\phi\colon M\to M'$
is an
	arrow in $\cat$. Then $\phi(\mathcal{N}) \smallsub M'$.
%\item \label{small_sum_1}
%		If $\mathcal{N} \smallsub M_1$ in $\cat$ and $M_2\in \cat$, then
%$i_1(\mathcal{N}) \smallsub M_1 \oplus M_2$.
\item 	\label{small_sum_2}
Let $\mathcal{N}_k\smallsub M$  in $\cat$, for $k$ in some finite index set $I$. Then
	$\sum_{k\in I} \mathcal{N}_k\smallsub M$.
\item \label{small_sum}
		Let $\mathcal{N}_k\subset M_k$,  $k \in I$, be a finite family of
		subobject-object pairs  in $\cat$. Then
		the following two assertions are equivalent:
		\begin{enumerate}
			\item  $\mathcal{N}_k \smallsub M_k$, for all $k\in I$;
			\item $\sum_{k\in I} i_k(\mathcal{N}_k) \smallsub
				\bigoplus_{k\in I}M_k$.
		\end{enumerate}
\end{enumerate}
\end{proposition}
\begin{proof}
\noindent\eqref{small_homomorphism}
Let $\mathcal{T} \subset M'$ be such that $\phi(\mathcal{N}) +
\mathcal{T} = \mathcal{M}'$.
Our first aim is to show that
 $\mathcal{N} +
\phi^{-1}(\mathcal{T}) = \mathcal{M}$.
Let $\psi_T \colon T \to M'$ and $\psi_N \colon N \to M$ be representatives of
$\mathcal{T}$ and $\mathcal{N}$, respectively. Denote by $\bcat$ the category
$\bcat(T,N,M,M')$, whose existence is asserted by Theorem~\ref{borceux}.
Since $\bcat$ is stable under finite limits and colimits, it is enough to show
that $\mathcal{N} + \phi^{-1}(\mathcal{T}) = \mathcal{M}$ in $\bcat$.

Let $R$ be a ring, such that there is a full and faithful exact functor
$F\colon \bcat \to R\Mod$.
Applying $F$ to $\phi(\mathcal{N}) + \mathcal{T} = \mathcal{M}'$ we get
 $(F\phi)(F\mathcal{N}) + F(\mathcal{T}) = F(M')$.
Now, for every $m \in F(M)$ there are $n \in F(\mathcal{N})$ and $t \in
F(\mathcal{T})$ such that $(F\phi)(m)= (F\phi)(n) + t$. In particular,
$(F\phi)(m -n ) = t $ belongs to $F(\mathcal{T})$. Therefore, $m -n \in
(F\phi)^{-1} ( F\mathcal{T})$. This shows that $m = n + (m-n)$ is an element of
$F(\mathcal{N}) + (F\phi)^{-1}(F \mathcal{T})$. Since $m$ was an arbitrary
element of $FM$, we get that $F(\mathcal{N}) + (F\phi)^{-1}(F\mathcal{T}) =
F(M)$. Since $F$ is exact this implies $\mathcal{N} +
\phi^{-1}(\mathcal{T}) = \mathcal{M}$.
As $\mathcal{N} \smallsub M$, we get $\phi^{-1}(\mathcal{T}) =
\mathcal{M}$. In particular, $\mathcal{N} \le \phi^{-1}(\mathcal{T})$.  Hence
$\phi (\mathcal{N}) \le \phi \phi^{-1}(\mathcal{T}) \le \mathcal{T}$. Therefore
$\mathcal{T} = \phi(\mathcal{N}) + \mathcal{T} = \mathcal{M}'$.
\medskip

\noindent \eqref{small_sum_2}  It is enough to prove the statement in  case the
cardinality of $ I $  is 2. Let $\mathcal{T} \subset M$ be such that
$\mathcal{N}_1 + \mathcal{N}_2 + \mathcal{T} = \mathcal{M}$.
	Since $\mathcal{N}_1 \smallsub M$, we get $\mathcal{N}_2 +
\mathcal{T} = \mathcal{M}$. Now $\mathcal{N}_2\smallsub M$,
	implies $\mathcal{T} =\mathcal{M}$.
\medskip

\noindent \eqref{small_sum}
	Suppose $\mathcal{N}_k\smallsub M_k$, for all $k\in I$. Then,
by~\eqref{small_homomorphism} we get
		$i_k (\mathcal{N}_k) \smallsub \bigoplus_{j\in I}
		M_j$. Therefore
$ \sum_{k \in I} i_k(\mathcal{N}_k) \smallsub
	\bigoplus_{k\in I} M_k$, by statement~\eqref{small_sum_2} of the
proposition.

Conversely, suppose $\sum_{k\in I} i_k(\mathcal{N}_k) \smallsub
\bigoplus_{k\in I} M_k$. Fix $\ell \in I$. We will show that $\mathcal{N}_\ell \smallsub
M_\ell$. Suppose $\mathcal{T} \subset M_\ell$ is such that $\mathcal{N}_\ell +
\mathcal{T} = \mathcal{M}_\ell$.

Define $\mathcal{S} = i_\ell(\mathcal{T}) + \sum_{k\not=\ell}
i_k (\mathcal{M}_k) $. Then $\sum_{k\in I} i_k(\mathcal{N}_k) +
\mathcal{S} = \sum_{k\in I} i_k(\mathcal{M}_k)$ is the top subobject of
$\bigoplus_{k\in I} M_k$.
Since $\sum_{k\in I} i_k(\mathcal{N}_k)$ is a superfluous subobject of
$\bigoplus_{k\in I} M_k$, we get $\mathcal{S} = \sum_{k\in I}
i_k(\mathcal{M}_k)$. Applying the $\ell$th canonical projection
$p_\ell \colon \bigoplus_{k\in I} M_k \to M_\ell$, we get
$\mathcal{T} = \mathcal{M}_\ell$. This shows that $\mathcal{N}_\ell \smallsub
M_\ell$.
\end{proof}
It should be noted that the properties stated in
Proposition~\ref{longlist}\eqref{small_sum_2} and \eqref{small_sum} cannot be
extended to infinite sums.
To give counter-examples we need the  notion of radical.
Given a ring $R$ and an $R$-module $M$ we can define $\rad(M)$ as the sum of all
superfluous subobjects in $M$. This definition is equivalent to the usual one
via the intersection of maximal subobjects by~\cite[Proposition~9.13]{fuller}.

Now consider the case $R=\Z$ and $M = \Q/\Z$. By~\cite[Exercise~9.2]{fuller},
we have $\rad(\Q/\Z) = \Q/\Z$. This shows that $\rad(\Q/\Z)$ is not a
superfluous submodule of $\Q/\Z$ despite being a sum of superfluous submodules.

Next let $R$ be a ring and $M$ an $R$-module such that $\rad(M)$ is not
superfluous in $M$. Denote by $S$ the set of all superfluous submodules in
$M$. The set $S$ is  infinite
as otherwise we would get a contradiction
to~Proposition~\ref{longlist}\eqref{small_sum}. Consider the submodule
$\widetilde{N} := \bigoplus_{N \in S} N$ of $\widetilde{M} :=
\bigoplus_{N \in S} M$. Then
$\widetilde{N}$ is  an infinite direct sum of superfluous submodules.
We will show that $\widetilde{N}$ is not superfluous.

Since $\rad(M)$ is not a superfluous subobject of $M$, there is a submodule $T\subset M$ such that $T\not=M$ and
$T + \rad(M) = M$.
We define
\begin{equation*}
\widetilde{T} = \left\{\, (m_N)_{N\in S} \in \widetilde{M} \,\middle|\,
\sum_{N \in S} m_N \in T \right\}.
\end{equation*}
We will show that $\widetilde{N} + \widetilde{T} = \widetilde{M}$ and
$\widetilde{T}\not=\widetilde{M}$. For the second assertion, note that if we take $m \in M
\setminus T$ and $N \in S$ then $i_N(m)$ is not an element of $\widetilde{T}$.

Now we will show that $\widetilde{T} + \widetilde{N} = \widetilde{M}$. For
this it is enough to check that for every $m \in M$ and $N\in S$ the element
$i_N(m)$ belongs to
$\widetilde{T} + \widetilde{N}$. Since $T +
\rad(M) = M$ and $\rad(M)$ is the sum of all superfluous subobjects we can write
$m$ as a linear combination
%\begin{equation*}
$m = t + \sum_{N'\in S'} m_{N'},\quad t \in T, \quad m_{N'} \in N'$,
%\end{equation*}
where $S'$ is a finite subset of $S$.
Now $i_{N'}(m_{N'}) \in \widetilde{N}$ for every $N'\in S'$. Further the  elements $\,i_N(t)\,$, $ \,i_N(m_{N'}) - i_{N'}(m_{N'})\,$ of $\widetilde{M}$
%\begin{equation*}
%\begin{aligned}
%i_N(t), \quad i_N(m_{N'}) - i_{N'}(m_{N'})
%\end{aligned}
%\end{equation*}
lie in $\widetilde{T}$. Therefore
\begin{equation*}
i_N(m) = \left( i_N(t) + \sum_{N'\in S'} (i_N(m_{N'}) - i_{N'}(m_{N'})) \right) +
\sum_{N' \in S'} i_{N'}(m_{N'})
\end{equation*}
is a sum of an element in $\widetilde{T}$ and of an element in
$\widetilde{N}$.
This shows $\widetilde{T} + \widetilde{N} = \widetilde{M}$ and thus
$\widetilde{N}$ is not a superfluous submodule of $\widetilde{M}$ despite
being
a direct sum of superfluous subobjects in the corresponding components.

For completeness of the exposition we recall the definitions of projective object and
projective resolution.
An object $P$ in $\cat$ is called \emph{projective} if for every epimorphism
$\phi\colon X\to Y$  in $\cat$ the map $\cat(P,f) \colon 	\cat(P,X)  \to
\cat(P,Y)$ is an epimorphism.
A \emph{projective cover} of $Y$ in $\cat$ is a projective object
$P$ together with an epimorphism $\pi\colon P\to Y$ such that $\ker\pi \smallsub
P$.
If a projective cover of $Y$ exists it is unique up to  isomorphism
(cf. \cite[Theorem~5.1]{semi-perfect}).

A \emph{projective resolution} of an object  $M\in \cat$ is an exact
complex $P_{\bullet} = (P_k, k\ge -1)$ with differentials $d_k \colon
P_{k+1} \to P_k$, $P_{-1} = M$, and $P_k$  projective.
This resolution  is called \emph{minimal} if, additionally,
 $\ker(d_k) \smallsub
P_{k+1}$ for all \mbox{$k\ge -1$}.

It follows from the definition, that if $P_\bullet$ is a minimal projective
resolution of $M$, then $d_{-1} \colon P_0\to M$ and $d_k\colon P_{k+1} \to
\ker(d_{k-1})$ for $k\ge 0$ are projective covers. Since a projective cover is
unique up to isomorphism, we see by an induction argument, that a minimal projective
resolution of $M$ is unique up to isomorphism.

\section{Graded algebras and modules}
\label{eilenberg}

 %In the case $\Gamma=\N$ these results
 %were originally proved in Eilenberg's

Let $R$ be a commutative ring with identity, $\Gamma$  a monoid with neutral element $\egamma$,
and $A$ a
$\Gamma$-\emph{graded associative $R$-algebra}, that is
we have an $R$-module decomposition
\begin{equation*}
	A = \bigoplus\limits_{\gamma\in \Gamma} A_{\gamma},
\end{equation*}
satisfying $A_{\gamma_1 }A_{\gamma_2 } \subset A_{\gamma_1 \gamma_2 }$ for all
$\gamma_1 $, $\gamma_2\in \Gamma $.
We also assume that the identity $e_A$ of $A$ is an element
of $A_{\egamma}$.

A left $A$-module $M$ is $\Gamma$-\emph{graded} if $M=
\bigoplus_{\gamma\in \Gamma} M_{\gamma}$, where each
$M_{\gamma}$ is an $R$-submodule of $M$, and  $
A_{\gamma_1} M_{\gamma_2} \subset M_{\gamma_1\gamma_2}$ for all $\gamma_1, \gamma_2 \in \Gamma.$
We will work with the category $A$-$\Gamma$-gr of left $\Gamma$-graded
$A$-modules and $A$-module homomorphisms respecting the
    grading: a map of $A$-modules $f\colon M_1\to M_2$ is in
$A$-$\Gamma$-gr if $f(M_{1,\gamma})\subset M_{2,\gamma}$ for all
$\gamma\in \Gamma$.
This category $A$-$\Gamma$-gr is abelian (see for example
Proposition~3.1 in \cite{semi-perfect}).

Given  a $\Gamma$-graded $A$-module $M$,  the support of $M$ is
\begin{equation*}
	\supp(M) = \left\{\, \gamma\in \Gamma \,\middle|\,  M_\gamma \not\cong
	0\right\}.
\end{equation*}

We call
a monoid $\Gamma$ equipped with an order $<$ an~\emph{ordered
monoid} if for any three elements $\alpha$, $\beta$, $\gamma\in
\Gamma$ such that $\alpha< \beta$ we have that also $\alpha\gamma
<\beta\gamma$ and $\gamma\alpha<\gamma\beta$.
Notice that we require that the multiplication with $\gamma$ preserves strict
inequalities. This is a stronger condition than to require just preservability
of non-strict inequalities. Terminological conventions on this matter can vary
from article to article.

A poset $(\mathcal{P},<) $  is called \emph{well-founded} if every strictly decreasing sequence of elements in
$\mathcal{P}$ is finite. With this definition we can speak about well-founded
ordered monoids.

The following theorem is a reformulation of Proposition~5.2 in \cite{perfect}.
\begin{theorem}
	\label{perfect}
	Let $\Gamma$ be a well-founded
ordered monoid such that
$\egamma$ is the least element of $\Gamma$.
	Then every $\Gamma$-graded $A$-module has a projective cover if and only if every $A_\egamma$-module has a projective cover.
\end{theorem}

Consider the general monoid $\Gamma$. For every $\beta\in \Gamma$ and every $\Gamma$-graded $A$-module $M$, we define the $\Gamma$-graded $A$-module
$M[\beta]$ to be $M$  as a left $A$-module with the homogeneous components
given
by
$\,M[\beta]_\gamma = \bigoplus_{\alpha \beta = \gamma} M_\alpha.
$
Note that each component $ M_\alpha$ of $M$ appears exactly once in the decomposition $M[\beta]=\bigoplus_{\gamma\in \Gamma} M[\beta]_\gamma$, since
 $M_\alpha$ is a
direct summand of $M[\beta]_{\alpha\beta}$ and not a direct summand of
any
$M[\beta]_{\gamma}$ for $\gamma\not= \alpha\beta$.

It is proved in~\cite[Proposition~4.3]{perfect}, that the set
$
\left\{\, A[\beta] \,\middle|\, \beta\in \Gamma \right\}
$
is a set of projective generators for the category of $\Gamma$-graded
$A$-modules. This fact will be used in the proof  of the following proposition.

\begin{proposition}
\label{ungrading_projective}
Let $\Gamma$ be a monoid and $A$ a $\Gamma$-graded algebra.
If $P$ is a projective $\Gamma$-graded $A$-module, then $P$ is projective as an $A$-module.
\end{proposition}
\begin{proof}
 Since $P$ is a projective $\Gamma$-graded $A$-module it is a direct summand of
$\bigoplus_{\beta \in \mathcal{B}} A[\beta]$, for some family $\mathcal{B}$ of
elements in  $\Gamma$. As $A[\beta] \cong A$ as an $A$-module we get that
$P$ is a direct summand of the free $A$-module $\bigoplus_{\beta\in
\mathcal{B}}A$. Hence $P$ is a projective $A$-module.
\end{proof}
Proposition~\ref{ungrading_projective} implies that if $P_\bullet$
is a projective resolution of a $\Gamma$-graded $A$-module $M$ then $P_\bullet$
is a projective resolution of $M$ considered as an $A$-module. If $P_\bullet$ is
minimal as $\Gamma$-graded projective resolution, it is not always true that it is minimal if  considered without grading.
Our next step will be to determine conditions under which minimality is
preserved upon forgetting the grading.

\begin{proposition}
	\label{ungrading_small}
	Let $\Gamma$ be a well-founded  ordered monoid with least element $\egamma$, and $A$ a $\Gamma$-graded algebra. Suppose
	$M$ is a $\Gamma$-graded $A$-module with finite support and $N$ is a
superfluous subobject of
	$M$ in the category $A$-$\Gamma$-gr. Then $N$ is a superfluous subobject of
	$M$ in the category $A$-Mod.
\end{proposition}
\begin{proof}
We will prove the proposition by induction on the cardinality of the support of
$M$. If $|
\supp(M)
|=1$ there is nothing to prove, as any $A$-submodule of $M$ is automatically $\Gamma$-graded.
Suppose the result holds for all $M$ with $|\sup(M)
|\le n-1$.
Consider $M$ with $|\supp(M)| = n$.
Let $T$ be an $A$-submodule of $M$ such that $N+T = M$. We have to show that
$T = M$.
Let $\gamma$ be a  maximal element of $\supp(M)$. Define $M'$ as the $\Gamma$-graded $A$-submodule of $M$ generated by
$\bigoplus_{\beta\not=\gamma} M_\beta$, that is $M'_\beta = M_\beta$ if
$\beta\not=\gamma$, and $M'_\gamma = \bigoplus_{\beta<\gamma} \bigoplus_{\alpha\beta
= \gamma} A_\alpha M_\beta$.
Note that in the last sum, we can take $\beta<\gamma$ and not just
$\beta\not=\gamma$, as the existence of $\alpha$ such that $\alpha \beta =
\gamma$ implies that
$\gamma > \beta$.
Also, since $\gamma$ is a maximal element of $\supp\left( M
\right)$,  both $M'_\gamma$ and $M_\gamma$ are $\Gamma$-graded
$A$-submodules of $M$. We will prove first that
\begin{equation}
	\label{first}
	M'_\gamma \subset T.
\end{equation}
By Proposition~\ref{longlist}\eqref{small_homomorphism}, we have that
\begin{equation}
	\label{nmgamma}
	\left.\raisebox{0.3ex}{$\left( N+M_\gamma
	\right)$}\middle/\!\raisebox{-0.3ex}{$M_\gamma$}\right.
	\smallsub
	\left.\raisebox{0.3ex}{$M$}\middle/\!\raisebox{-0.3ex}{$M_\gamma$}\right.
\end{equation}
in the category $A$-$\Gamma$-gr. Since the support of
$\left.\raisebox{0.3ex}{$M$}\middle/\!\raisebox{-0.3ex}{$M_\gamma$}\right. $ has
cardinality $n-1$, by the induction hypothesis we get that
\eqref{nmgamma} also holds in $A$-Mod. Obviously
\begin{equation*}
	\left.\raisebox{0.3ex}{$\left(N+M_\gamma\right)$}\middle/\!\raisebox{-0.3ex}{$M_\gamma$}\right.
	+
	\left.\raisebox{0.3ex}{$\left( T+M_\gamma
	\right)$}\middle/\!\raisebox{-0.3ex}{$M_\gamma$}\right.
	=
	\left.\raisebox{0.3ex}{$M$}\middle/\!\raisebox{-0.3ex}{$M_\gamma$}\right..
\end{equation*}
Therefore
\begin{equation}
	\label{m1prime}
	\left.\raisebox{0.3ex}{$\left( T+M_\gamma
	\right)$}\middle/\!\raisebox{-0.3ex}{$M_\gamma$}\right.
	=
	\left.\raisebox{0.3ex}{$M$}\middle/\!\raisebox{-0.3ex}{$M_\gamma$}\right..
\end{equation}
It is now easy to show \eqref{first}.
In fact, we only have to check that $A_\alpha M_\beta\subset T$
for every $\beta<\gamma$ and $\alpha$, such that $\alpha\beta = \gamma$.
Note that, for such $\alpha$ we have $\alpha\not= \egamma$ and so $\alpha>\egamma$.
Let $y\in M_\beta$  and $a \in A_\alpha$. It
follows from \eqref{m1prime} that there is $z \in M_\gamma$ such that $y + z \in
T$. Thus also $a(y+z) \in T$. Since $\alpha \gamma >\gamma$ and $\gamma$ is a
maximal element of $\supp(M)$, we get that $az =0$. Thus $ay \in T$,
i.e. $A_\alpha M_\beta \subset T$.
Now define
\begin{align*}
	\overline{ M }& =
	\left.\raisebox{0.3ex}{$M$}\middle/\!\raisebox{-0.3ex}{$M'_\gamma$}\right.,
	& \overline{ N } &=
	\left.\raisebox{0.3ex}{$\left( N+M'_\gamma
	\right)$}\middle/\!\raisebox{-0.3ex}{$M'_\gamma$}\right.
	\subset \overline{ M }
	, &\overline{ T } &=
	\left.\raisebox{0.3ex}{$T$}\middle/\!\raisebox{-0.3ex}{$M'_\gamma$}\right.
	\subset \overline{ M }
	.
\end{align*}
Note that $\overline{M}$ is the
internal direct sum of $\overline{M}_\gamma$ and $\bigoplus_{\beta\not=\gamma}
\overline{M}_\beta$ in the category of $\Gamma$-graded $A$-modules.
Our next step is to prove that $\overline{ N } \smallsub \overline{ M }$ in the
category $A$-Mod.
 By Proposition~\ref{longlist}\eqref{small_homomorphism} we have that $\overline{ N }\smallsub \overline{
M }$ in the category $A$-$\Gamma$-gr.
Since $\overline{N}$ is $\Gamma$-graded, we get that
$\overline{N}_\gamma$ is a $\Gamma$-graded $A$-submodule of $\overline{M}_\gamma$ and
$\bigoplus_{\beta\not=\gamma} \overline{N}_\beta$ is a $\Gamma$-graded $A$-submodule
of $\bigoplus_{\beta\not=\gamma} \overline{M}_\beta$. Moreover $ (\bigoplus_{\beta\not=\gamma} \overline{N}_\beta) \bigoplus \overline{N}_\gamma = \overline{N}$.
Therefore from Proposition~\ref{longlist}\eqref{small_sum}, we have that $\overline{N}_\gamma \smallsub
\overline{M}_\gamma$ and $\bigoplus_{\beta\not=\gamma} \overline{N}_\beta \smallsub
\bigoplus_{\beta\not=\gamma} \overline{M}_\beta$ in the category $A$-$\Gamma$-gr.
Since the cardinalities of the supports of
$\overline{M}_\gamma$ and
of $\bigoplus_{\beta\not=\gamma} \overline{M}_\beta$ are less than $n$,
the induction assumption gives that
$\overline{N}_\gamma \smallsub
\overline{M}_\gamma$ and $\bigoplus_{\beta\not=\gamma} \overline{N}_\beta \smallsub
\bigoplus_{\beta\not=\gamma} \overline{M}_\beta$ in
$A$-Mod. Applying Proposition~\ref{longlist}\eqref{small_sum} with $\cat$ the category $A$-Mod, we conclude
that $\overline{N} \smallsub \overline{M}$ in $A$-Mod. Thus
$
	\overline{N} + \overline{ T } = \overline{ M}
$
implies that
$ \overline{ T } = \overline{M}$. Therefore $T = M$.
\end{proof}
We will use the result of Proposition~\ref{ungrading_small} only in the case
$M$ is a projective module. It is natural to wonder if the condition of finiteness on
$|\supp(M)|$ is redundant. The next example shows
that this is not the case.

Let $R$ be a commutative ring with identity, with the property that  its Jacobson radical $J=J(R)$ is not left
T-nilpotent. For example, we can take $R = \Z_{(p)}$,  the localization of $\Z$ at the ideal $(p)$, for some prime $p$.   Consider $R$ as an
$\N$-graded ring with $R_0 = R$ and all other homogeneous components equal to zero. Define the
$\N$-graded $R$-module $F$ by $F_n = {}_R R$. Clearly $F$ is a projective
$\N$-graded $R$-module.
If we consider $F$ as an $R$-module, by Proposition~9.19 in~\cite{fuller}, we have
 $N := \rad(F)=\bigoplus_{n \in \N} J$.
In particular, $N$ is an $\N$-graded submodule of $F$. It is easy to check that
$N \smallsub F$ in $R$-$\N$-gr. Indeed, if $T = \bigoplus_{n \in \N}T_n$ is another
$\N$-graded $R$-submodule of $F$ such that $N + T  = F$, then for each
component $n$ we have $N_n + T_n = R$, i.e. $J + T_n =R$. Since $J\smallsub R$ in $R$-Mod, this implies $T_n =R$ and $T = F$.

By Proposition~17.10 and Lemma~28.3 in~\cite{fuller}, we know that $N =  JF$, and that $JF$ is a superfluous subobject of $F$ in $R$-Mod if and
only if $J$ is a left $T$-nilpotent ideal. Hence, under our assumption on the ring
$R$, $N$ is not a superfluous submodule of $F$ in $R$-Mod. This shows that
forgetting the grading can render a superfluous subobject to become
non-superfluous.

\begin{proposition}
\label{finsupres}
	Let $\Gamma$ be a monoid,
	$A$ a $\Gamma$-graded algebra with finite support, and $M$ a
$\Gamma$-graded $A$-module with finite support.
Then there exists a projective resolution
 $P_\bullet$ of $M$ in $A$-$\Gamma$-gr such that each $P_k$ has  finite
support.
\end{proposition}
\begin{proof}
It is enough to show that for every
$M$ with finite support there is a projective $\Gamma$-graded $A$-module
$P$ with finite support and an epimorphism $f\colon P \to M$. Then using this fact for $\ker(f)$, we get $P_2 \to \ker (f)$, and hence the first two steps of a projective resolution of $M$ with finite support.
Repeating this process recursively we obtain a projective resolution of the required type.

Since $\left\{\, A[\beta] \,\middle|\, \beta\in\Gamma \right\}$ is a set of
projective generators of $A$-$\Gamma$-gr, there is $P= \bigoplus_{\beta\in \Gamma}
A[\beta]^{\kappa_\beta}$, where $\kappa_\beta$ are cardinals, and an
epimorphism of $\Gamma$-graded $A$-modules $f\colon P\to M $.
Suppose $\beta\not\in\supp(M)$. Then the restriction of $f$ to each summand $A[\beta]$
in $P$ is zero. In fact, the module $A[\beta]$ is generated as an  $A$-module by
the element $e\in A_e \subset A[\beta]_\beta$, and  its image under
$f$ is in $M_\beta=0$.
Therefore, without loss of generality, we can assume that $\kappa_\beta=0$ for
all $\beta\not\in \supp(M)$.
Since $|\supp(A[\beta])| \leq |\supp(A)|$  and $A$ has finite support, we get that $|\supp(P)|\leq  \sum_{\beta\in
\supp(M)} |\supp(A[\beta])|  \le |\supp(M)| \cdot | \supp(A) |$ is finite.
\end{proof}

We can now state and prove the main result of this section.

\begin{theorem}
	\label{minimal_ungrading}
	Let $\Gamma$ be a well-founded ordered monoid with least element $\egamma$, and
	$A$ a $\Gamma$-graded algebra with finite support. Given a
$\Gamma$-graded $A$-module $M$ with finite support, let $P_\bullet$
	be a minimal projective resolution of
	$M$ in $A$-$\Gamma$-gr. Then $P_\bullet$ is a minimal projective
	resolution of $M$ in the category $A$-Mod.
In other words, the grading forgetting functor from $A$-$\Gamma$-gr to $A$-Mod
preserves minimal projective resolutions of $\Gamma$-graded $A$-modules with finite support.
\end{theorem}
\begin{proof} We know that $P_\bullet$ is a projective resolution of $M$ in the category
$A$-Mod by Proposition~\ref{ungrading_projective}. Thus we have only to check
that it is minimal.

By Proposition~\ref{finsupres} there is a projective resolution $\overline{ P
}_\bullet$ of $M$ in
$A$-$\Gamma$-gr such that all $\overline{ P }_k$, $k\ge 0$, have finite
support.
Since $P_\bullet$ is a minimal projective resolution (by applying for example Theorem~5.1
in \cite{semi-perfect}) there is an embedding of $P_\bullet$ into $\overline{
P }_\bullet$.
 Thus each
$P_k$, $k\ge 0$, has finite support.
 Since the resolution $P_\bullet$ is  minimal in $A$-$\Gamma$-gr, all the maps
$d_{-1}\colon P_0\to M$  and $d_k\colon P_{k+1} \to P_k$ for $k\ge 0$ have
superfluous kernels  in $A$-$\Gamma$-gr. From
Proposition~\ref{ungrading_small}, these kernels are also superfluous in $A$-Mod.
\end{proof}

\section{Twisted products}
We start this section with an overview of the concept of twisted product of
rings. Then we specialise to twisted products of $\Gamma$-graded
algebras and modules, and study under which conditions  the functor $- \ltimes_{\Gamma} N$, defined below, preserves minimal projective resolutions.

Given rings $S$, $A_1$, and $A_2$ , suppose we have   ring homomorphisms $\phi_i\colon S\to A_i$, for
$i=1,2$. We say that $A$ is a \emph{twisted product}  of $A_1$ and
$A_2$ over $S$ if there are a ring homomorphism $\phi \colon S \to A$ and  an
$S$-bimodule isomorphism  $\gamma\colon A_1 \otimes_S A_2 \to A$ such that
\begin{equation}
	\label{axioms}
	\begin{aligned}
	\gamma(\phi_1(s)\otimes 1) &= \gamma(1\otimes \phi_2(s)) = \phi(s),&
\gamma (a_1\otimes a_2) & = \gamma(a_1\otimes 1) \gamma(1\otimes a_2) \\[2ex]
\gamma(a_1 a_1' \otimes 1) &= \gamma(a_1\otimes 1) \gamma(a_1'\otimes
	1),
	& \gamma(1\otimes a_2 a_2' ) &= \gamma(1\otimes a_2) \gamma(1\otimes
	a_2').
\end{aligned}
\end{equation}
If $A$ is a twisted product of $A_1$ and $A_2$ over $S$, one can define a
\emph{twisting
homomorphism} of abelian groups
\begin{equation*}
\begin{aligned}
	T\colon A_2 \otimes_S A_1 & \to A_1\otimes_S A_2\\
	a_2\otimes a_1 & \mapsto \gamma^{-1} (\gamma(1\otimes a_2)
	\gamma(a_1\otimes 1)).
\end{aligned}
\end{equation*}
Note that it is then possible to reconstruct  $A$ from $\phi_1$, $\phi_2$ and
$T$. The name~\emph{twisted product} is justified by the existence of the map $T$.

Twisted products of algebras over fields where studied in~\cite{cap}.
A more general approach, that can be applied to monoids in arbitrary monoidal
categories, was considered in~\cite{beck}.

Next we study a twisted product involving a $\Gamma$-graded
algebra. As usual,  all the unnamed tensor products are considered over $R$.

We say that an $R$-algebra $B$ is a $\Gamma$-\emph{algebra} if there is a
right action of $\Gamma$ on $B$
\begin{align*}
 r\colon B\times \Gamma&\to B\\
 (b,\gamma)  & \mapsto b^{\gamma},
\end{align*}
 such that for each
$\gamma\in \Gamma$ the map $b\mapsto b^\gamma$
is an algebra homomorphism.

Let $A$ be a $\Gamma$-graded $R$-algebra and $B$ a $\Gamma$-algebra.
We define a binary operation $m$ on $A\otimes B$ by
\begin{equation*}
(a \otimes b)(a' \otimes b')=a a'\otimes b^\beta b', \mbox{ for }  a \in
A_{\alpha},\  a' \in A_{\beta},\  b, b' \in B.
\end{equation*}
The $R$-module $A\otimes B$ when considered together with the binary operation $m$ will  be denoted by
$A\ltimes_{\Gamma} B$. It is routine to check that the following proposition holds.

\begin{proposition}
The pair $(A\ltimes_{\Gamma} B,m)$ is an $R$-algebra with identity $1_A
\otimes 1_B$. It is a twisted product of
$A$ and $B$ (over R), where $\phi_A$, $\phi_B$,  and $\phi_{A\ltimes_\Gamma B}$ are the
unity maps, and $\gamma\colon A\otimes B \to A\ltimes_\Gamma B$ is the identity map.
Moreover $A\ltimes_{\Gamma}B$ is
 $\Gamma$-graded, with  the grading given by
$\left(A\ltimes_{\Gamma} B\right)_\gamma =A_{\gamma}\otimes B, \, \, \gamma\in
\Gamma.$
\end{proposition}
Note that the maps $A\to A\ltimes_{\Gamma} B$ and $B\to
A\ltimes_{ \Gamma} B$  given by
\begin{align*}
a & \mapsto a\otimes 1_B &  b & \mapsto 1_A \otimes b
\end{align*}
	are  homomorphisms of algebras, being the first one a homomorphism of $\Gamma$-graded algebras.

Let $N$ be a $B$-module and $M= \bigoplus_{\gamma\in
\Gamma}M_{\gamma}$ a $\Gamma$-graded $A$-module. We
define an $(A\ltimes_{\Gamma}B)$-module structure on $M\otimes N$ as
follows
$$
a_{\gamma_1}\otimes b \otimes m_{\gamma_2}\otimes x \mapsto
a_{\gamma_1}m_{\gamma_2}\otimes b^{\gamma_2}x,
$$
for all $a_{\gamma_1} \in A_{\gamma_1}$, $b\in B$, $m_{\gamma_2}\in
M_{\gamma_2}$ and $x\in N$. We denote this module by
$M\ltimes_{\Gamma}N$.
This is a $\Gamma$-graded module, with $(M\ltimes_{\Gamma}
N)_{\gamma} = M_{\gamma}\otimes N$.

 Let $\varphi\colon M_1\to M_2$ be a homomorphism of $\Gamma$-graded $A$-modules
 and $\psi\colon N_1\to N_2$ a homomorphism of $B$-modules. We write
 $\varphi\ltimes_{\Gamma}\psi$  for the map
 \begin{align*}
  M_1\ltimes_{\Gamma}N_1& \to M_2\ltimes_{\Gamma}N_2\\
  m\otimes x  & \mapsto \varphi(m)\otimes \psi(x).
 \end{align*}
 %\begin{proposition}
Clearly $\varphi\ltimes_{\Gamma} \psi$ is a homomorphism of
  $\Gamma$-graded $A\ltimes_{\Gamma} B$-modules.
It follows  that the correspondence
\begin{align*}
 (M,N)& \mapsto M\ltimes_{\Gamma}N\\
 (\varphi,\psi)& \mapsto \varphi\ltimes_{\Gamma} \psi
\end{align*}
gives a bifunctor from the categories $A$-$\Gamma$-gr and $B$-Mod to the
category $(A\ltimes_{\Gamma}B)$-$\Gamma$-gr.
In particular, for each $B$-module $N$, we have the functor $-\ltimes_{\Gamma}N$
from the category $A$-$\Gamma$-gr to the category
$(A\ltimes_{\Gamma}B)$-$\Gamma$-gr. This functor is exact if and only if
$N$ is a flat $R$-module. Note that
$-\ltimes_\Gamma N$ preserves arbitrary direct sums.
In fact, the $R$-isomorphism
\begin{equation*}
\begin{aligned}
\phi \colon \Big(\bigoplus_{i\in I} M_i \Big) \ltimes_\Gamma N & \to \bigoplus_{i\in I}
(M_i \ltimes_\Gamma N) \\
(m_i)_{i\in I} \otimes x & \mapsto (m_i \otimes x)_{i\in I}.
\end{aligned}
\end{equation*} is also  an isomorphism of
$\Gamma$-graded
$A\ltimes_\Gamma B$-modules.

Next we will establish a sufficient condition for the functor $- \ltimes_{\Gamma} N$ to
preserve projective objects.
For every $\beta \in \Gamma$, we denote by ${}_\beta N$ the $B$-module with the
same underlying abelian group as $N$ but with the $B$-action defined by
\begin{equation}
\label{dodot}
b \cdot_\beta x = b^\beta x,
\end{equation}
for every $x\in N$ and $b\in B$.

\begin{proposition}
\label{projective}
Let $A$ be a $\Gamma$-graded $R$-algebra and $B$ a $\Gamma$-algebra. Suppose that
$N$ is a $B$-module with the property that all the $B$-modules ${}_\beta N$, $\beta \in
\Gamma$, are projective. Then, for any   projective $\Gamma$-graded $A$-module $P$, the $\Gamma$-graded $A\ltimes_{\Gamma} B$-module $P \ltimes_\Gamma
N$ is   projective.
\end{proposition}
\begin{proof}
First we consider the case $P= A\left[ \beta \right]$, for some
$\beta\in\Gamma$. We claim that $A\left[ \beta \right] \ltimes_\Gamma N\cong
\left( A\ltimes_\Gamma {}_\beta N \right)\left[ \beta \right]$ as
$\Gamma$-graded $A\ltimes_\Gamma B$-modules, and so it is projective. For this note that, for any $a'\in A_\alpha$,  $a\in A$, $b\in B$, and $x\in N$, the $A\ltimes_{\Gamma} B$-action  on $A\left[ \beta \right] \ltimes_\Gamma N$  gives $(a \otimes b) (a'\otimes x)  =  aa' \otimes
b^{\alpha\beta}  x.$ On the other hand, in $\left( A\ltimes_\Gamma {}_\beta N \right)\left[ \beta \right]$ there holds $(a \otimes b) (a'\otimes x)  =  aa' \otimes
(b^{\alpha})^\beta  x=aa' \otimes
b^{\alpha\beta}  x.$ Therefore the identity map gives the desired isomorphism.
 %$\phi$ from $A[\beta]\ltimes_{\Gamma} N $ to $\left(

Now let $P$ be an arbitrary
projective $\Gamma$-graded $A$-module. Then $P$  is a direct summand of
$\bigoplus_{\beta \in I} A[\beta]$, for  some family $I$ of
elements in $\Gamma$.
Since
 the functor $-\ltimes_\Gamma N$  preserves  direct sums, we get that
$P\ltimes_\Gamma N$ is a direct summand of
\begin{equation*}
\begin{aligned}
\Big( \bigoplus_{\beta\in I} A\left[ \beta \right] \Big)\ltimes_\Gamma N
\cong
\bigoplus_{\beta\in I}\Big(  A\left[ \beta \right] \ltimes_\Gamma N\Big),
\end{aligned}
\end{equation*}
which is projective. Therefore $P\ltimes_\Gamma N$ is a
projective $\Gamma$-graded $A\ltimes_\Gamma B$-module.
\end{proof}

Let $N$ be a $B$-module which is flat as an $R$-module and such that all ${}_\beta N$
are projective $B$-modules. Then Proposition~\ref{projective} shows
that the functor $-\ltimes_{\Gamma}N$ preserves projective resolutions.
Note that it does not map in general a minimal projective resolution
into a minimal projective resolution.

\begin{proposition}
	\label{small_N}
	Suppose  $M_1 \smallsub M_2$ in
	$A$-$\Gamma$-gr.
Then,  for any $B$-module N which is finitely generated over $R$, we have   $M_1 \ltimes_\Gamma N \smallsub M_2 \ltimes_\Gamma N$ in $A\ltimes_\Gamma
B$-$\Gamma$-gr.
\end{proposition}
\begin{proof}
	Let $T$ be a $\Gamma$-graded $A\ltimes_\Gamma B$-submodule of
	$M_2\ltimes_\Gamma N$ such that
	\begin{equation}
		\label{Y}
		M_1\ltimes_\Gamma N + T = M_2 \ltimes_\Gamma N.
	\end{equation}
Every $\Gamma$-graded $A\ltimes_\Gamma B$-module can be considered as a $\Gamma$-graded $A$-module via the
canonical homomorphism $A\to A\ltimes_{\Gamma} B$.
	Therefore
	\eqref{Y}
	also holds in the category of $\Gamma$-graded $A$-modules. Let
	$\left\{ x_1, \dots, x_k \right\}$ be a generating set of $N$ over
$R$.  Since $M_1 \smallsub M_2$ in
	$A$-$\Gamma$-gr, we get from Proposition~\ref{longlist}\eqref{small_sum}, that
	$\bigoplus_{j=1}^k M_1 \otimes R x_j \smallsub \bigoplus_{j=1}^k M_2
	\otimes R x_j$ in $A$-$\Gamma$-gr.
Consider the canonical epimorphism
\begin{equation*}
\begin{aligned}
\phi \colon \bigoplus_{j=1}^k M_2 \otimes R x_j \to M_2\ltimes_\Gamma N
\end{aligned}
\end{equation*}
 of $\Gamma$-graded $A$-modules. We have  $\phi\left( \bigoplus_{j=1}^k M_1 \otimes Rx_j
\right) = M_1 \ltimes_\Gamma N$.
	Thus, by Proposition~\ref{longlist}\eqref{small_homomorphism},
$M_1 \ltimes_\Gamma N \smallsub M_2 \ltimes_{\Gamma} N$ in
$A$-$\Gamma$-gr.  Therefore $T = M_2\ltimes_\Gamma N$.
\end{proof}

The following corollary is an immediate consequence of Propositions~\ref{projective} and~\ref{small_N}.
\begin{corollary}
	\label{minimalmap}
	Let  $\phi \colon P\to M$ be
	projective cover in $A$-$\Gamma$-gr. Suppose that  $N$  is a  $B$-module which is flat
and finitely generated over $R$, and such that all ${}_\beta N$ are projective
$B$-modules.  Then $\phi\ltimes_\Gamma N\colon P\ltimes_\Gamma
	N
	\to M\ltimes_\Gamma N$ is a projective cover in
	$A\ltimes_\Gamma B$-$\Gamma$-gr.
\end{corollary}
For future reference, we bring together in the next theorem the results proved in this section.
\begin{theorem}\label{to be stated}
Let $A$ be a $\Gamma$-graded algebra and $B$ a $\Gamma$-algebra. Suppose that $N$ is a
$B$-module which is flat and finitely generated over $R$, and  such that for all $\beta \in \Gamma$ the $B$-modules ${}_\beta N$
are projective. Then the functor \mbox{\rm $-\ltimes_\Gamma N\colon A\ggr \to
A\ltimes_\Gamma B\ggr $} preserves minimal projective resolutions.
\end{theorem}

\section{Relative stratifying ideals and projective resolutions}
\label{RSIPR}
In this section we adopt a different setting.
Once more $R$ denotes a commutative ring with identity, but $A$ is simply an associative $R-$algebra.
Given an ideal $I$ of $A$, we are interested in determining conditions for the  functor $A/I \otimes_A - \colon  A\md \rightarrow  A/I\md$ (or, equivalently, the functor $N\mapsto N/IN $)
to preserve minimal projective resolutions.
 For this we will use relative homological algebra. So we start with a brief overview of this topic.

 We say that an $A$-module
$P$ is
$(A,R)$-\emph{projective} if
for every epimorphism $f\colon M\to N$ which is split as an epimorphism of
$R$-modules, the homomorphism $\Hom_A(P,f)$ is surjective.

Given an $A$-module $M$ and an exact complex $P_\bullet \onto M$, we say that
$P_\bullet$ is an $(A,R)$-\emph{projective resolution of} $M$ if every $P_k$ is an
$(A,R)$-projective module and the complex $P_\bullet\onto M$ is split as a
complex of $R$-modules.
Every $A$-module $M$ admits a canonical $(A,R)$-projective resolution,
$
\beta_\bullet\left( A,R,M \right)\onto M,
$
known as \emph{bar} resolution. Recall that
$
\beta_k\left( A,R,M \right) = A^{\otimes_R (k+1)}\otimes_R M,
$
and the differentials and the splitting maps are the usual ones and can be found
in~\cite{homology}.

Let $N$ be a right $A$-module and $M$ a left $A$-module. Given an
$(A,R)$-projective resolution $P_\bullet$ of $M$, we define the \emph{relative tor groups}
$
\Tor_k^{(A,R)}(N,M) = H_k(N\otimes_A P_\bullet)
$, all $k\ge 0$.
It follows from Theorem~IX.8.5 in~\cite{homology}, that the groups
$\Tor_k^{(A,R)}(N,M)$ are independent of the choice of the $(A,R)$-projective
resolution of $M$ and, in particular, can be computed using the bar resolution
of $M$.

As we mentioned in the introduction, in~\cite{LU} we obtained an efficient
combinatorial criterion for a triple $(A,R,I)$ to have the property
$\Tor^{(A,R)}_k (A/I,A/I) \cong 0$ for $k\ge 1$. In the next series of
propositions
we derive various consequences of this property under additional conditions on
$(A,R,I)$ culminating in Theorem~\ref{cor:unrelative}. This gives a criterion
for the functor $A/I\otimes -$ to preserve (minimal) projective
resolutions.

\begin{proposition}
\label{projcond}
Let $A$ be an $R$-algebra  and  $I$ an
 ideal of $A$ such that, for $k\ge 1$,  $\Tor^{(A,R)}_k (A/I,A/I) \cong 0$.
Suppose that $A$ and $A/I$ are projective as right $R$-modules. Then  $\Tor^{(A,R)}_k (A/I,M) \cong 0$, for~$k\ge1$ and  any $M \in A/I\md $.
\end{proposition}
\begin{proof}
Since $\Tor^{(A,R)}_k (A/I,A/I) \cong 0$ for $k\ge 1$, we get that the complex
\begin{equation*}
A/I \otimes_A \beta(A,R,A/I)\onto A/I\otimes_A A/I
\end{equation*}
is exact. Moreover, the differentials in this
complex are homomorphisms of $A/I$-bimodules.  Further, the first term of this complex is
$A/I \otimes_A A/I \cong A/I$ and  every other term  is
of the form
$
A/I \otimes_A A^{\otimes (k+1)}\otimes A/I \cong A/I \otimes A^{\otimes k}
\otimes A/I,
$
where all the unnamed tensor products are over $R$ and $k\ge 0$.
Now, since $A/I$ is a projective right $R$-module, we get that $A/I \otimes A$ is a projective
right $A$-module. This fact together with the fact that  $A$ is a projective right $R$-module, implies that $A/I \otimes
A$ is a projective right $R$-module. Continuing, we get that $A/I\otimes
A^{\otimes k} $ are projective right $R$-modules, for all $k\ge 0$. Thus $A/I \otimes
A^{\otimes k} \otimes A/I$ is a right $A/I$-projective  module.
Therefore the exact complex $A/I\otimes_A \beta(A,R,A/I)\onto A/I\otimes_A A/I $ % is exact, all differentials
splits in the category of right
$A/I$-modules. Hence
$
A/I\otimes_A \beta(A,R,A/I)\otimes_{A/I} M \onto A/I\otimes_A A/I
\otimes_{A/I} M
$
is an exact
complex. But it is isomorphic to
$
A/I \otimes_A \beta(A,R,M) \onto A/I\otimes_A M \cong M,
$
and therefore it
computes the torsion groups $\Tor^{(A,R)}_k (A/I,M)$. We get then
 that $\Tor^{(A,R)}_k (A/I,M) \cong
0, \, k\ge 1$.
\end{proof}
In the next proposition we relate relative  with classical torsion groups.
\begin{proposition}
\label{relab}
Let $A$ be a free $R$-algebra, $I$ an
ideal of $A$, and $M  $ an $R$-free left $A$-module. Then $\Tor^{(A,R)}_k(A/I,M)
\cong \Tor^A_k(A/I,M)$, for all
$k$.
\end{proposition}
\begin{proof}
We consider the bar resolution $\beta(A,R,M)$ of $M$. Every
module in this resolution is of the form $A\otimes A^{\otimes k} \otimes M$,
with $k\ge 0$, where all
the tensor products are taken over $R$. Since $M$ and $A$ are free
$R$-modules,  $A^{\otimes k} \otimes  M$ is a free
 $R$-module. Hence $A\otimes A^{\otimes k} \otimes  M $ is a free $A$-module. This shows that
$\beta(A,R,M)$ is a projective resolution of $M$ in the category of left
$A$-modules.  Now, both tor groups $\Tor^{(A,R)}_k (A/I,M)$ and $\Tor^A_k
(A/I,M)$ can be computed using the complex $A/I\otimes_A
\beta(A,R,M)$. This proves the result.
\end{proof}
\begin{theorem}
\label{cor:unrelative} Let  $M \in  A/I\md$ be an $R$-free left module.
Then, in the conditions of the previous two propositions, the functor  $A/I \otimes_A -$
sends every projective resolution of $ M$  in $A\md$  to a projective
resolution of $M$ in $A/I\md$. If the initial resolution  of $M$  in $A\md$ is minimal, then the final resolution  in $A/I\md$ is also a minimal projective resolution of $M$ .
\end{theorem}
\begin{proof}  Let $P_\bullet\onto M$ be a projective resolution of
$M$ in $A$-Mod.
By Propositions~\ref{projcond} and~\ref{relab},  $\Tor^A_k(A/I,M)
\cong 0$, for $k\ge 1$. Therefore, since
 $\Tor^A_0 (A/I,M) = A/I \otimes_{A} M \cong M$,
 the complex $A/I \otimes_A P_\bullet \onto M$ is exact. As every
$P_k$ is a projective $A$-module, it follows that $A/I \otimes_A P_k$ is an
$A/I$-projective module.

 To prove that the minimality is preserved,  consider in $A$-Mod the minimal projective resolution $P_\bullet\onto M$, with differentials  $d_k\colon P_{k+1} \to
P_k$ for $k\ge -1$ (for simplicity, we write $P_{-1}=M)$. This can be decomposed into  short exact sequences
$$ 0 \rightarrow \ker \alpha_k \xrightarrow{\beta_k} P_{k+1} \xrightarrow{\alpha_k} \ker \alpha_{k-1} \rightarrow 0,$$ where $ d_k=\beta_{k-1}\alpha_k. $ Write $F$ for the functor $A/I \otimes_A - $. Then it is easy to see that $\ker (F(d_k)) = \ker( F(\alpha_k))=  \image (F(\beta_k))= \pi  (\image (\beta_k))=  \pi  (\ker (\alpha_k))=  \pi  (\ker (d_k)),$ where $\pi\colon  P_{k+1} \mapsto F(P_{k+1})$ is the epimorphism given by $x \mapsto 1\otimes x$. Since $ \ker (d_k)\smallsub
P_{k+1}$, the result follows from  Proposition~\ref{longlist}~\eqref{small_homomorphism}.

\end{proof}

%right adjoint + linear  implies preserves projectives?
\section*{Acknowledgment}
The authors are grateful to Leonid Positselski for the
discussion~\cite{posic} which led to the example after
Proposition~\ref{ungrading_small}.

\section*{Funding}
This work was partially supported by the Centre for Mathematics of the
University of Coimbra
- UIDB/00324/2020, funded by the Portuguese Government through FCT/MCTES.
The second author was also partially  supported via the FCT grant
CEECIND/04092/2017.
\bibliography{kostant}
\bibliographystyle{amsplain}
\end{document}